\documentclass[12pt]{amsart}
\usepackage{amsmath,amssymb}
\usepackage{mathrsfs}
\usepackage{verbatim}
\usepackage{enumerate}
\usepackage{multicol}
\usepackage{mathtools}
\usepackage{tikz}
\usepackage{fancyhdr}
\usepackage{enumerate}

\newcommand{\bigchi}{\raisebox{2pt}{\large$\chi$}}

\newtheorem{thm}{Theorem}

\newtheorem{prop}[thm]{Proposition}

\newtheorem{lemma}[thm]{Lemma}

\theoremstyle{definition}

\newtheorem{remark}[thm]{Remark}
\newcommand{\abar}{\mathbf{a}}
\newcommand{\xbar}{\mathbf{x}}
\newcommand{\bbar}{\mathbf{b}}
\newcommand{\ebar}{\mathbf{e}}
\newcommand{\cbar}{\mathbf{c}}

\newcommand{\spn}{{\text{span}}}
\newcommand{\field}{\Bbbk}
\newcommand{\maxideal}{\mathfrak{m}}

\newcommand{\diff}{\partial}

\newcommand{\sym}{\mathfrak S}
\newcommand{\symn}{\mathfrak S_n}

\newcommand{\cosetspace}{V}

\newcommand{\ann}{\mathrm{Ann}}

\DeclareMathOperator{\gr}{gr}
\DeclareMathOperator{\LF}{LF}

\DeclareMathOperator{\soc}{Soc}
\DeclareMathOperator{\rk}{rk}

\newcommand{\apol}{S}
\DeclareMathOperator{\image}{im}
\DeclareMathOperator{\rank}{rank}

\newcommand{\hf}{\mathrm{HF}}
\newcommand{\hs}{\mathrm{HS}}

\begin{document}

\begin{abstract}
We describe the graded characters and Hilbert functions of certain graded artinian Gorenstein quotients of the polynomial ring which are also representations of the symmetric group. Specifically, we look at those algebras whose socles are trivial representations and whose principal apolar submodules are generated by the sum of the orbit of a power of a linear form.
\end{abstract}

\title[Hilbert Functions of $\symn$-Stable Artinian Gorenstein Algebras]{Hilbert Functions of\\ $\symn$-Stable Artinian Gorenstein Algebras}

\author[Geramita]{Anthony V. Geramita}
\address{
Department of Math \& Stats\\
Jeffery Hall, University Ave.\\
Kingston, ON Canada, K7L 3N6\newline\indent
Dipartimento di Matematica\\
Via Dodecaneso, 35\\
16146 Genova, Italy
}
\email{tony@mast.queensu.ca}

\author[Hoefel]{Andrew H. Hoefel}
\address{
Department of Math \& Stats\\
Jeffery Hall, University Ave.\\
Kingston, ON Canada, K7L 3N6
}
\email{ahhoefel@mast.queensu.ca}

\author[Wehlau]{David L. Wehlau}
\address{Department of Mathematics and Computer Science \\
Royal Military College \\ Kingston, Ontario, Canada \\ K7K 5L0
}
\email{wehlau@rmc.ca}



\maketitle
\section{Introduction} \label{sec:intro}
Let $R = \field[x_1, \ldots, x_n] = \bigoplus_{k \geq 0} R_k$ be the standard graded polynomial ring in $n$ variables over a field $\field$ (which is one of $\mathbb Q$, $\mathbb R$, or $\mathbb C$)
and let $\symn$ denote the symmetric group on $n$ letters. We are interested in the Hilbert functions and graded characters of graded artinian Gorenstein algebras which are also representations of $\symn$. Specifically, we will examine quotients of $R$ whose one-dimensional socles are spanned by a symmetric polynomial $F$.

Every homogeneous polynomial $F$ of degree $d$ can be expressed as a linear combination  of $d$-th powers of linear forms $L_1,\ldots,L_m$ for some $m$. If we symmetrize this expression for $F$ by summing over all permutations and dividing by $n!$, we can express $F$ as a linear combination of $F_1,\ldots,F_m$ where each $F_i$ is the sum of the elements in the orbit of $L_i^d$. In this paper, we will only consider polynomials $F$ which are the sum of the orbit of the $d$-th power of a single linear form (i.e., $m=1$).

 For example, let $n=4$ and consider the linear form
 \[ L=x_1 + x_2 + 2x_3 + 3x_4 \ .
\]
The degree $7$ symmetric polynomial
 \[ F = \sigma_1 L^7 + \cdots + \sigma_{12} L^7  \]
spans the one-dimensional socle of the graded artinian Gorenstein algebra $R/I_F$ where $I_F$ consists of all $f \in R$ with $\diff F/ \diff f = 0$.

The dimensions of the homogeneous components of $R/I_F$ are recorded in its Hilbert function $\hf_{R/I_F}(k) := \dim_\field (R/I_F)_k$, although it is often convenient to present them in a generating function called a Hilbert series:
\[ \hs_{R/I_F}(t) := \sum_{k \geq 0} \dim_\field (R/I_F)_k t^k.
\]

Currently, there is no known description of all Hilbert functions that arise from graded Gorenstein algebras (except for $n=2,3$ \cite[Theorem 1.44]{MR1735271} \cite{MR0485835}).
In order to find the Hilbert function of an artinian graded algebra with a given socle, one typically computes the ranks of a collection of Catalecticant matrices --- one for each degree \cite{MR1702103}.

In our example, however, $R/I_F$ has additional structure which
will allow us to find that its Hilbert series is
\[ \hs_{R/I_F}(t) = 1 + 4t + 9t^2 +  12t^3 + 12t^4 + 9t^5 + 4t^6 + t^7.
\]

Our point of view is illustrated as follows: since $F$ is a symmetric polynomial, $I_F$ is stable under the action of $\sym_n$. Thus, $I_F$ and its homogeneous components $(I_F)_k$ are representations of $\symn$. So, the quotient $R/I_F = \bigoplus_{k \geq 0} (R/I_F)_k$, where $(R/I_F)_k = R_k/(I_F)_k$, is a graded representation of $\symn$. If we denote the \emph{character} of a finite dimensional representation $V$ of $\symn$ by $\bigchi_V: \symn \to \field$,
then the \emph{graded character} of $R/I_F$ is defined by
\[ \bigchi_{R/I_F}(t) := \sum_{k\geq0} \bigchi_{(R/I_F)_k} t^k.
\]
This encodes the algebra's structure as a graded representation much like the Hilbert series does for its structure as a graded vector space.

Recall that since any representation of $\symn$ is a direct sum of irreducible representations, and the irreducible representations of $\symn$ are in one-to-one correspondence with partitions $\lambda \vdash n$, $\lambda = (\lambda_1,\ldots, \lambda_r)$, $\lambda_1 \geq \cdots \geq \lambda_r \geq 1$, we can write
\[ \bigchi_{(R/I_F)_k} = \sum_{\lambda \vdash n} m_\lambda \bigchi^\lambda
\]
where $m_\lambda \in \mathbb N$ and $\bigchi^\lambda$ is the character of the irreducible respresentation corresponding to $\lambda$.


Writing $\bigchi^{(\lambda_1,\ldots,\lambda_r)}$ as $\bigchi^{\lambda_1\cdots\lambda_r}$,
the graded character of $R/I_F$, for $F$ as above, is
\begin{align*}
	\bigchi_{R/I_F}(t)
	= \qquad &\bigchi^4 \\
	{}+ (&\bigchi^4 + \bigchi^{31})t \\
	{}+ (&\bigchi^4 +  2\bigchi^{31} + \bigchi^{22})t^2 \\
	{}+ (&\bigchi^4 +  2\bigchi^{31} + \bigchi^{22} + \bigchi^{211})t^3\\
	{}+ (&\bigchi^4 +  2\bigchi^{31} + \bigchi^{22} + \bigchi^{211})t^4\\
	{}+ (&\bigchi^4 +  2\bigchi^{31} + \bigchi^{22})t^5 \\
	{}+ (&\bigchi^4 + \bigchi^{31})t^6 \\
	{}+\phantom{(} &\bigchi^4 t^7.
\end{align*}

In this paper, we describe the graded characters of Gorenstein algebras $R/I_F$ whose socles are spanned by a form $F$ which is the sum of the orbit of a power of a linear form $L$. (Our only additional requirements are that the coefficients of $L$ are real and that they do not sum to zero. The latter requirement is essentially equivalent to the embedding dimension of $R/I_F$ being $n$.) The graded characters of such algebras are palindromic, as in the example above, and the multiplicities of the irreducible characters in each degree can be computed directly from certain Hall-Littlewood polynomials, which is the content of our main theorem (see Theorem \ref{thm:maintheorem}).

  In fact, the graded character of $R/I_F$ depends on the degree of $F$ and the number of repeated coefficients in $L$, but not on the values of those coefficients. Furthermore, the Hilbert function of $R/I_F$ can be recovered from the graded character by replacing each $\bigchi^\lambda$ occurring in $\bigchi_{R/I_F}(t)$ with the numbers $f^\lambda$ which count the number of standard Young tableaux with shape $\lambda$ and which is equal to the dimension of the irreducible representation of type $\lambda$. For an introduction to the representation theory of the symmetric group we recommend \cite{MR1824028} (and \cite{MR2538310} to see how these representations are realized in the polynomial ring).

Since artinian Gorenstein algebras are characterized by having a one-dimensional socle,
it follows that graded artinian Gorenstein algebras which admit an action of $\symn$ come in two types, which correspond to the only two one-dimensional representations of $\symn$. The socle of the algebra is either
\begin{enumerate}[(i)]
\item the trivial representation and spanned by a symmetric polynomial, or
\item the alternating representation and spanned by an alternating polynomial.
\end{enumerate}

In \cite{MR1779775}, Bergeron, Garsia and Tesler
described the graded character of artinian Gorenstein algebras of type (ii).
The graded character of such an algebra is a multiple of the graded character of the coinvariant algebra $R_\mu$, where $\mu = (1,\ldots, 1) \vdash n$. Roth \cite{MR2134302} extended their result to any artinian algebra whose socle is spanned by alternating forms.

Morita, Wachi and Watanabe \cite{MR2542134} found the Hilbert function of each isotypic piece of the algebra
$A(n,k) = \field[x_1,\ldots,x_n]/(x_1^k,\ldots,x_n^k)$ for each $n$ and $k$. The algebra $A(n,k)$ is of type $(i)$ since its socle is spanned by the symmetric monomial $F = (x_1\cdots x_n)^{k-1}$.

We will consider algebras $R/I_F$ of type (i)
where we have chosen the socle of $R/I_F$ to be spanned by a symmetric polynomial $F$ of a specific type.

The paper is organized as follows: in Section 2 we review some facts about Gorenstein graded algebras, the action of $\symn$ on $R$ and the apolarity module and its submodules.

In Section \ref{section:orbit sum} we apply the ideas of Section 2 to a specific type of $\symn$-stable submodules of $R$. 

In Section \ref{section:point} we demonstrate our main theorem, which shows that the graded characters of the algebras we are considering are determined by certain algebras that were considered by, among others, DeConcini, Garsia and Procesi.  This connection allows us to calculate the graded characters of our algebras.

In the final section we suggest some possibilities for further inquiry.

\textbf{Acknowledgments.} The authors would like to thank Henry de Valence and Dave Stringer who computed the first examples of graded characters for this project as part of an NSERC summer undergraduate research assistantship. We also thank Federico Galetto for many helpful conversations concerning $R_\mu$ and Zach Teitler for his helpful comments on the first version of this paper.

All three authors gratefully acknowledge the partial support of NSERC for this work.

\section{Gorenstein Algebras and Apolarity Modules}  The following elementary remarks are well known and can be found, for example, in the expository article \cite{MR1381732}.

A standard graded algebra $A = R/I$ is \emph{Gorenstein} if it contains a maximal $A$-regular sequence which generates an irreducible ideal of $A$.

If $A$ is \emph{artinian} (i.e., finite dimensional as a $\field$ vector space) then the socle of $A$ is the ideal $\soc A = ( 0 : \maxideal) = \{ f \in A \mid f \maxideal = 0 \}$ where $\maxideal = \bigoplus_{k \geq 1} A_k$ is the homogeneous maximal ideal of $A$.  If $A$ is artinian and Gorenstein then its socle must be one-dimensional as a vector space.

Let each $\sigma \in \symn$ act on a linear form $a_1x_1 + \cdots + a_nx_n \in R_1$ by
\begin{align*}
	\sigma(a_1 x_1 + \cdots + a_n x_n)
	&= a_1 x_{\sigma(1)} + \cdots + a_n x_{\sigma(n)} \\
	&= a_{\sigma^{-1}(1) }x_1 + \cdots + a_{\sigma^{-1}(n)}x_n.
\end{align*}
This action extends to an action on $R$ which is given by
\[ (\sigma f)(x_1, \ldots, x_n) = f(x_{\sigma(1)}, \ldots, x_{\sigma(n)}).
\]

If we write a monomial in $R$ as $\xbar^\bbar = x_1^{b_1} \cdots x_n^{b_n}$, for an exponent vector $\bbar = (b_1, \ldots, b_n) \in \mathbb N^n$, then
\begin{align*}
\sigma(\xbar^{\bbar})
&= x_{\sigma(1)}^{b_1} \cdots x_{\sigma(n)}^{b_n}=x_1^{b_{\sigma^{-1}(1)}} \cdots x_n^{b_{\sigma^{-1}(n)}}.
\end{align*}
Accordingly, we define $\sigma(\bbar) = (b_{\sigma^{-1}(1)}, \ldots, b_{\sigma^{-1}(n)})$ for $\sigma \in \symn$ and
for any $\bbar \in \mathbb N^n$ or $\field^n$, so that $\sigma \xbar^\bbar = \xbar^{\sigma \bbar}$.
One can check that all of the above are left actions.

Moreover, if $\abar \in \field^n$ then
$$f(\sigma \abar) = f(a_{\sigma^{-1}(1)}, \ldots, a_{\sigma^{-1}(n)}) = (\sigma^{-1}f)(\abar) \ .
$$
So, if $\ebar_i$ is the $i$-th standard basis vector of $\field^n$, then $\sigma(\ebar_i) = \ebar_{\sigma(i)}$.

Partial differentiation is a $\field$-bilinear operator $\diff: R \times R \to R$ where
the partial derivative of the monomial $\xbar^\cbar$ by the monomial $\xbar^\bbar$ is defined to be
\begin{align*}
	\diff(\xbar^\bbar, \xbar^\cbar)
	&= \frac{c_1! \cdots c_n!}{(c_1-b_1)! \cdots (c_n - b_n)!} \xbar^{\cbar - \bbar} \\
	&:= \frac{\cbar!}{(\cbar - \bbar)!} \xbar^{\cbar - \bbar}
\end{align*}
when $\cbar \geq \bbar$ (i.e., $c_i \geq b_i$ for all $i$) and zero otherwise.

We  extend this definition linearly in both components to define $\diff(f,g)$, the partial derivative of a polynomial $g$ by another polynomial $f$.

Partial differentiation endows $R$ with an $R$-module structure which is different from that of a rank-one free module. To avoid confusion,
we let $\apol=R$ be this $R$-module where for  $f \in R$ and $g \in \apol$, we set $fg := \diff(f,g) \in \apol$.  $\apol$ is called the \emph{apolarity module} of $R$. So the bilinear operator $\diff$ can be viewed as a map $\diff: R \times S \to S$.

Write  $R = \bigoplus_{k\geq 0} R_k$, using the standard grading on $R$.    Let $\apol = \bigoplus_{k \geq 0} \apol_k$   where $\apol_k = R_k$ as vector spaces. We consider $\apol$ as a graded $R$-module even though its graded components satisfy the slightly unconventional condition that $R_k \apol_j \subseteq \apol_{j-k}$ for all $j,k \in \mathbb N$ where $\apol_{j-k} = 0$ for $k > j$. This convention allows a polynomial $f$ to have the same degree regardless of whether it is in $R$ or $\apol$.

The bilinear operator $\diff$ restricts to a bilinear form $\diff_k: R_k \times \apol_k \to \field$ which  induces two maps $R_k \to \apol_k^*$ and $\apol_k \to R_k^*$. Let $\overline f$ denote the polynomial obtained by taking the complex conjugate of the coefficients of the polynomial $f$. Since $\diff_k(f,\overline f) \neq 0$ for any $f \in R_k \setminus \{0\}$, the map $R_k \to \apol_k^*$ is injective and hence is an isomorphism (and similarly for the map $\apol_k \to R_k^*$). Thus, the bilinear form $\diff_k$ is a perfect pairing.

For any subspace $M_k$ of $\apol_k$, we define
\[ M_k^\perp = \{ f \in R_k \mid \forall g \in M_k,\, \diff_k(f, g) = 0 \}.
\]
Then $\diff_k$ induces a well-defined bilinear form
\[ \diff_{k,M_k}: R_k/M_k^\perp \times M_k \to \field
\]
which is also a perfect pairing. Thus, $M_k^*$ and $R_k/M_k^\perp$ are isomorphic vector spaces and hence $M_k$ and $R_k/M_k^\perp$ have the the same dimension.

If we now consider a graded submodule $M \subseteq \apol$, its annihilator is the homogeneous ideal
\[ \ann(M) = \{ f \in R \mid \forall g \in M,\, \diff(f,g) = 0  \}
\]
and so $R/\ann(M)$ is a graded algebra.  The homogeneous components of $\ann(M)$ are $\ann(M)_k = M_k^\perp$ (see \cite[Proposition 2.5]{MR1381732}).  Thus the Hilbert function of $M$ (i.e., $\hf_M(k) = \dim_\field M_k$) is equal to the Hilbert function of $R/\ann(M)$.

Just as submodules $M$ of $S$ determine quotients $R/\ann(M)$ of $R$ by taking $\ann(M)_k = M_k^\perp$ for $k \geq 0$, homogeneous ideals of $R/\ann(M)$ determine graded quotients of $M$. In particular, the socle of $R/\ann(M)$ determines $M/\maxideal M$ (any basis of which represents a set of minimal generators of $M$) since the homogeneous components of $\soc(R/\ann(M))$ and $\maxideal M$ are orthogonal under $\diff_{k,M_k}$.

This fact, proved by Macaulay \cite[\S60]{MR1281612} (and referred to as the Macaulay Correspondence), shows that artinian Gorenstein algebras are in one-to-one correspondence with principally generated submodules of the apolarity module. (See \cite[Lemma 2.12]{MR1735271} for a modern treatment of these facts.)

We would now like to consider submodules of $\apol$ which are also graded representations of $\symn$.
We have already defined an action of $\symn$ on $R$ and, using the same action, $\apol$ is also a representation of $\symn$. Since the homogeneous components of $R$ and $\apol$ are stable under the action of $\symn$, $R$ and $\apol$ are graded representations.
Partial differentiation is an invariant bilinear form in that $\diff_k(\sigma f, \sigma g) = \diff_k(f,g)$ for all $f \in R_k$, $g \in \apol_k$ and $\sigma \in \symn$. Therefore, the dual representation of $M_k$ is equivalent to $(R/\ann(M))_k$. (The action of $\symn$ on the dual $V^*$ of a representation $V$ is given by $(\sigma f)(v) = f(\sigma^{-1} v)$ for $\sigma \in \symn$, $f \in V^*$ and $v \in V$.) Since $\sigma$ and $\sigma^{-1}$ are conjugate, the symmetric group has the special property that all of its representations are self-dual. Thus, $R/\ann(M)$ and $M$ have the same graded characters.


Suppose $M = \langle g \rangle$ is the principal submodule of $\apol$ generated by $g$. As a vector space, $M$ consists of all partial derivatives of $g$. In this case, $R/\ann(M)$ is an artinian Gorenstein algebra which is isomorphic to $M$ as an $R$-module. If $g$ is a homogeneous polynomial of degree $d$, then $M$ is a graded $R$-module and its homogeneous components are given by $M_k = \image \theta_{g,k}$ where $\theta_{g,k}: R_{d-k} \to \apol_k$ is the map
\[ \theta_{g,k}(f) = \diff(f,g).
\]
We set $\theta_g: R \to \apol$ to be the map $\theta_g = \bigoplus_{k \geq 0} \theta_{g,k}$ which can also be described by $\theta_g(f) = \diff(f,g)$. This map $\theta_g$ is called the \emph{Catalecticant map} of $g$ and its image is $M = \bigoplus_{k \geq 0} M_k = \langle g \rangle$.

\section{The Apolarity Module of an Orbit Sum}\label{section:orbit sum}

Let $M = \langle g \rangle$ be a principally generated submodule of  $\apol$  where $g$ is a {\it symmetric} homogeneous polynomial $g$ of degree $d$.   Then $M$ is an $\symn$-stable subspace of $\apol$ and
\text{$\theta_{g,k}: R_{d-k} \to \apol_k$}
is equivariant. Thus, $M_k$ and $R_{d-k}/\ker \theta_{g,k} = (R/\ann(M))_{d-k}$ are equivalent representations.

As $M_{d-k}$ is also equivalent to $(R/\ann(M))_{d-k}$, the graded character of $M$ is palindromic in the sense that $\bigchi_{M_k} = \bigchi_{M_{d-k}}$ for all $k$.



We now restrict our attention to a specific family of symmetric functions.  Let $L$ be a linear form and $d>0$ a positive integer.  Let $\{ \sigma_1 L, \dots , \sigma_\ell L \}$ be the orbit of $L$. Define
$$
F= \Sigma _{i =1}^\ell \sigma_i(L^d) .
$$
We will consider the artinian Gorenstein ring $R/I_F$ where $I_F$ is the annihilator of $M = (F) \subset S$.

Suppose $\abar = (a_1, \ldots, a_n) \in \field^n$ and let
\[ L = a_1 x_1 + \cdots + a_n x_n \in \apol_1.
\]
Let $b_1, \ldots, b_r$ be all the distinct coordinates of $\abar$.
Define $$\mu_j = | \{ i \mid a_i = b_j\} |$$ for $1 \leq j \leq r$.
Reorder $b_1,\ldots,b_r$ so that $\mu_1 \geq \mu_2 \geq  \cdots  \geq \mu_r \geq 1$ and let
$\mu = (\mu_1, \ldots, \mu_r)$. Thus the partition $\mu \vdash n$ \emph{associated to} $\abar$ characterizes the number of repeated coordinates in $\abar$.

As mentioned above, $\symn$ acts on the left of tuples by $\sigma(\abar) = (a_{\sigma^{-1}(1)},\ldots, a_{\sigma^{-1}(n)})$. The stabilizer of $\abar$ is
\[ (\symn)_\abar = \{ \sigma \in \symn \mid \sigma \abar = \abar\}.
\]
Let $\ell$ be the index of the stabilizer of $\abar$ in $\symn$.  Then  $\ell = \binom{n}{\mu} := n!/(\mu_1!\cdots \mu_r!)$.

If $\sigma_1, \dots ,\sigma_\ell$ are left coset representatives of the stabilizer of $\abar$ in $\symn$ then let $F \in \apol_d$ be the homogeneous symmetric polynomial determined by $\abar$, i.e.
\begin{equation}
F =  \sum_{i=1}^\ell \sigma_i L^d.
\label{eq:Fdef}
\end{equation}




In order to obtain the Hilbert function of $R/I_F$ explicitly, where $F$ is as in (\ref{eq:Fdef}) above, we use the fact that $\hf_{R/I_F} = \hf_{M}$ and instead examine
\[ M = \bigoplus_{k \geq 0} \image \theta_{F,k}
\]
where   $\theta_{F,k}(f) = \diff(f,F)$.

\begin{lemma}\label{lem:thetaformula}
For $F = \sum_{i=1}^\ell \sigma_i L^d \in \apol_d$ as described above and $\xbar^\bbar$, a monomial of degree $d-k$, we have
\[
	\theta_{F,k}(\xbar^\bbar)
	= \frac{d!}{k!} \sum_{i=1}^\ell \abar^{\sigma_i^{-1}(\bbar)} \sigma_i L^k.
\]
\begin{proof} The proof is by induction on $d-k$.
If $d-k=0$ then $\bbar = (0,\ldots,0)$ and the formula obviously holds.

Assume that the formula above holds for every monomial $\xbar^\bbar$ of some fixed degree $d-k$. Differentiating our formula by $x_j$ gives
\begin{align*}
	\diff(\xbar^{\bbar + \ebar_j},F)
	&= \diff(x_j,\, \diff(\xbar^\bbar, F))\\
	&= \diff\left(x_j,\ \frac{d!}{k!} \sum_{i=1}^\ell \abar^{\sigma_i^{-1}(\bbar)} \sigma_i L^k\right) \\
	&= \frac{d!}{k!} \sum_{i=1}^\ell \abar^{\sigma_i^{-1}(\bbar)} \diff(x_j, \sigma_i L^k)\\
	&= \frac{d!}{k!} \sum_{i=1}^\ell \abar^{\sigma_i^{-1}(\bbar)} k (\sigma_i L^{k-1}) \diff(x_j, \sigma_i L).
\end{align*}
Since $\sigma_i L = a_{\sigma_i^{-1}(1)}x_1 + \cdots + a_{\sigma_i^{-1}(n)}x_n$, we have $\diff(x_j, \sigma_i L) = a_{\sigma_i^{-1}(j)}$. Furthermore,
$a_{\sigma_i^{-1}(j)} = \abar^{\ebar_{\sigma_i^{-1}(j)}} = \abar^{\sigma_i^{-1}(\ebar_j)}$ where $\ebar_j$ is the $j$-th standard basis vector of $\field^n$.
So, $\abar^{\sigma_i^{-1}(\bbar)} a_{\sigma_i^{-1}(j)} = \abar^{\sigma_i^{-1}(\bbar)} \abar^{\sigma_i^{-1}(\ebar_j)} = \abar^{\sigma_i^{-1}(\bbar + \ebar_j)}$.
Thus,
\begin{align*}
	\diff(\xbar^{\bbar + \ebar_j},F)
	&= \frac{d!}{k!} \sum_{i=1}^\ell \abar^{\sigma_i^{-1}(\bbar)} k (\sigma_i L^{k-1}) \diff(x_j, \sigma_i L)\\
	&= \frac{d!}{(k-1)!} \sum_{i=1}^\ell \abar^{\sigma_i^{-1}(\bbar + \ebar_j)} \sigma_i L^{k-1},
\end{align*}
proving the formula holds for any monomial $\xbar^{\bbar+\ebar_j}$ of degree ${d- k+1}$.
\end{proof}
\end{lemma}

Recall that $\sigma_1,\ldots, \sigma_\ell$ are representatives for the left-cosets of $(\symn)_\abar$.
We define the $\ell$-dimensional vector space $\cosetspace$ as the span of these representatives of the left-cosets of the stabilizer of $\abar$. Then $\cosetspace$ is a representation of $\symn$ where $\tau(\sigma_i) := \sigma_j$ if $\tau \sigma_i \in \sigma_j (\symn)_\abar$.

For any given degree $k$, we define $\phi_k : \cosetspace \to \apol_k$ by setting
\[  \phi_k(\sigma_i) = \sigma_i L^k
\]
for $1 \leq i \leq \ell$ and extending linearly. The definition of $\phi_k$ does not depend on our choice of coset representatives since if $\tau \in \sigma_i (\symn)_\abar$, then $\tau = \sigma_i \gamma$ for some $\gamma \in (\symn)_\abar$ and hence $\tau L^k = \sigma_i \gamma L^k = \sigma_i L^k$. For similar reasons, $\phi_k$ is equivariant.

Let $\psi_k : R_k \to \cosetspace$ be the linear map given by
\[ \psi_k(\xbar^\bbar) = \sum_{i=1}^\ell \binom{k}{\bbar} \abar^{\sigma_i^{-1}\bbar} \sigma_i.
\]

\begin{lemma} \label{lem:rankpsi} For all $k \geq 0$, $\rank \psi_k = \rank \phi_k$.
\begin{proof}
The coefficient of $\xbar^\bbar \in \apol_k$ appearing in $\phi_k(\sigma_i)$ is equal to the coefficient of $\sigma_i^{-1} \xbar^\bbar = \xbar^{\sigma_i^{-1}\bbar}$ in $L^k$, which is $\binom{k}{\bbar} \abar^{\sigma_i^{-1}\bbar}$.
If we order the monomial bases of $R_k$ and $\apol_k$ in the same manner, then the matrices of $\phi_k$ and $\psi_k$ are transposes of each other, and thus have the same rank.
\end{proof}
\end{lemma}

\begin{lemma} The map $\psi_k$ is equivariant.
\begin{proof}
For $\tau \in \symn$,
$\tau \psi_k(\xbar^\bbar) = \sum_{i=1}^\ell \binom{k}{\bbar} \abar^{\sigma_i^{-1} \bbar} \tau \sigma_i$. In this expression, $\tau \sigma_i$ holds the place of the left coset $\tau \sigma_i (\symn)_\abar$.
Let $\sigma_j$ be our chosen representative for this coset. Thus, $\tau \sigma_i (\symn)_\abar = \sigma_j (\symn)_\abar$ and hence $\tau \sigma_i = \sigma_j \gamma$ for some $\gamma \in (\symn)_\abar$.
Rearranging we get $\sigma_i^{-1} = \gamma^{-1} \sigma_j^{-1} \tau$
and therefore $\abar^{\sigma_i^{-1}\bbar} = \abar^{\sigma_j^{-1} \tau \bbar}$, since $\abar^{\gamma \cbar} = \abar^\cbar$ for any $\gamma$ in the stabilizer of $\abar$ and any exponent vector $\cbar$.
Thus,
\[ \tau \psi_k(\xbar^\bbar) = \sum_{i=1}^\ell \binom{k}{\bbar} \abar^{\sigma_i^{-1} \bbar} \tau \sigma_i = \sum_{j=1}^\ell \abar^{\sigma_j^{-1} \tau \bbar}\sigma_j = \psi_k(\tau \bbar)
\]
\end{proof}
\end{lemma}

Finally, let $\nu_k: R_k \to R_k$ be the non-singular equivariant linear scaling map defined by
\[ \nu_k(\xbar^\bbar) = \binom{k}{\bbar}^{-1} \xbar^\bbar= \frac{b_1! \cdots b_n!}{k!} \xbar^\bbar.
\]

We thus have the following sequence of maps:
\[ R_{d-k} \stackrel{\nu_{d-k}}{\longrightarrow} R_{d-k} \stackrel{\psi_{d-k}}{\longrightarrow} V \stackrel{\phi_k}{\longrightarrow} \apol_k.
\]

The following theorem gives the ``factored'' presentation of $\theta_{F,k}$ promised earlier.

\begin{thm} \label{thm:theta} Let $F = \sum_{i=1}^\ell \sigma_i L^d \in \apol_d$ and let $\theta_{F,k}: R_{d-k} \to \apol_k$,  $\nu_{k}$, $\phi_k$ and $\psi_k$ be as given above. For all $k$ with $0 \leq k \leq d$, we have
\[	\theta_{F,k} = \frac{d!}{k!}\, \phi_k \circ \psi_{d-k} \circ \nu_{d-k}.
\]
\begin{proof}
Applying the composition $\phi_k \circ \psi_{d-k} \circ \nu_{d-k}$ to a monomial $\xbar^\bbar \in R_{d-k}$ and using linearity gives
\begin{align*}
\phi_k \circ \psi_{d-k} \circ \nu_{d-k}(\xbar^\bbar)
&= \binom{d-k}{\bbar}^{-1} \phi_k \circ \psi_{d-k}(\xbar^\bbar)\\
&= \binom{d-k}{\bbar}^{-1} \sum_{i=1}^\ell \binom{d-k}{\bbar} \abar^{\sigma_i^{-1}\bbar}\phi_k(\sigma_i)\\
&= \sum_{i=1}^\ell \abar^{\sigma_i^{-1}\bbar} \sigma_i L^k.
\end{align*}
The result then follows from Lemma \ref{lem:thetaformula}.
\end{proof}
\end{thm}




The Hilbert function of $R/I_F$ is equal to the Hilbert function of the image of $\theta_F$. Thus, to determine $\hf_{R/I_F}(k)$ it suffices to find the rank of $\theta_{F,k} = \frac{d!}{k!} \phi_k \circ \psi_{d-k} \circ \nu_{d-k}$. Since $\nu_{d-k}$ is non-singular, it suffices to determine the rank of the composition of $\phi_k$ and $\psi_{d-k}$. To that end, we need to examine the relationship between $\image \psi_{d-k}$ and $\ker \phi_k$.

Using the distinguished basis $\sigma_1,\ldots, \sigma_\ell$ of $\cosetspace$, we can define the dot
product of two vectors $v = c_1 \sigma_1 + \cdots + c_\ell \sigma_\ell$ and $w = d_1 \sigma_1 + \cdots + d_\ell \sigma_\ell$ to be $v\cdot w = c_1 \overline{d_1} + \cdots + c_\ell \overline{d_\ell}$.

In the following results, when we require $\abar$ to have real coordinates, we still allow $\field$ to be one of $\mathbb Q$, $\mathbb R$, or $\mathbb C$.

\begin{prop}\label{prop:perpkerimg} If the coordinates of $\abar$ are real then $\ker \phi_k$ and $\image \psi_k$ are orthogonal complements with respect to the dot product on $V$. In particular, their intersection is trivial.
\begin{proof}
Take $v \in \ker \phi_k$ and $w \in \image \psi_k$. If $v = \sum_{i=1}^\ell c_i \sigma_i$, then $\phi_k(v) = \sum_{i=1}^\ell c_i \sigma_i L^k = 0$. The coefficient of $\xbar^\bbar$ in $\sum_{i=1}^\ell c_i \sigma_i L^k$ is $\sum_{i=1}^\ell c_i \binom{k}{\bbar} \abar^{\sigma_i^{-1} \bbar}$ and this coefficient must be zero for each monomial $\xbar^\bbar$ of degree $k$.

Since $w \in \image \psi_k$, there is a homogeneous polynomial
$f = \sum_{\bbar} r_\bbar \xbar^{\bbar}$ of degree $k$ with $w = \psi_k(f) = \sum_\bbar r_\bbar \sum_{i=1}^\ell \binom{k}{\bbar} \abar^{\sigma_i^{-1} \bbar} \sigma_i$. Thus, $w = \sum_{i=1}^\ell d_i \sigma_i$ where $d_i = \sum_\bbar r_\bbar \binom{k}{\bbar} \abar^{\sigma_i^{-1} \bbar}$. Since $\abar$ is real, $\overline{d_i} = \sum_\bbar \overline{r_\bbar} \binom{k}{\bbar} \abar^{\sigma_i^{-1} \bbar}$.
Therefore,
\begin{align*}
	v \cdot w
	&= \sum_{i=1}^\ell c_i \overline{d_i} \\
	&= \sum_{i=1}^\ell c_i \sum_\bbar \overline{r_\bbar} \binom{k}{\bbar} \abar^{\sigma_i^{-1} \bbar}\\
	&= \sum_\bbar \overline{r_\bbar} \sum_{i=1}^\ell c_i \binom{k}{\bbar} \abar^{\sigma_i^{-1} \bbar}\\
	&= 0,
\end{align*}
since $\sum_{i=1}^\ell c_i \binom{k}{\bbar} \abar^{\sigma_i^{-1} \bbar} = 0$ for each $\xbar^\bbar$ of degree $k$. Furthermore, $\phi_k$ and $\psi_k$ have the same rank since their matrices are transposes of each other. Thus, $(\ker \phi_k)^\perp = \image \psi_k$.
\end{proof}
\end{prop}

\begin{lemma} \label{lemma:kercontainment} If the coordinates of $\abar$ are real and $a_1 + \cdots + a_n \neq 0$ then $\ker \phi_{k+1} \subseteq \ker \phi_k$ and $\image \psi_k \subseteq \image \psi_{k+1}$ for all $k \geq 0$.
\begin{proof}
Using Proposition \ref{prop:perpkerimg}, it suffices to show the containment of the kernels.
If we suppose $v = \sum_{i=1}^\ell c_i \sigma_i \in \ker \phi_{k+1}$, then $\sum_{i=1}^\ell c_i \sigma_i L^{k+1} =0$. Therefore $\diff(x_1+\cdots + x_n,\, \phi_{k+1}(v))$ is both zero and
\begin{align*}
	\diff(x_1+\cdots + x_n,\, \phi_{k+1}(v))
	&= \sum_{i=1}^\ell c_i \diff(x_1+ \cdots +x_n,\, \sigma_i L^{k+1})\\
	&= \sum_{i=1}^\ell c_i (k+1) (a_1 + \cdots + a_n)  \sigma_i L^k\\
	&= (a_1 + \cdots + a_n) (k+1) \phi_k(v).
\end{align*}
Since we have assumed that $a_1+ \cdots + a_n \neq 0$, we have $v \in \ker \phi_k$. Thus, $\ker \phi_{k+1} \subseteq \ker \phi_k$.
\end{proof}
\end{lemma}

\begin{prop} If the coordinates of $\abar$ are real and $a_1 + \cdots + a_n \neq 0$ then for any $i, j \geq 0$, $\ker \phi_i$ and $\image \psi_j$ meet transversely.
\begin{proof}
Since $\ker \phi_j$ and $\image \psi_j$ are orthogonal complements, we have $\ker \phi_j + \image \psi_j = V$ and $\ker \phi_j \cap \image \psi_j = 0$.
If $i \geq j$ then, by Lemma \ref{lemma:kercontainment}, $\ker \phi_i \subseteq \ker \phi_j$ and hence $\ker \phi_i \cap \image \psi_j = 0$.
If $i < j$ then, by Lemma \ref{lemma:kercontainment}, $\ker \phi_j \subseteq \ker \phi_i$ and hence $\ker \phi_i + \image \psi_j = V$.
\end{proof}
\end{prop}

\begin{prop} \label{prop:gradedcharM}
If the coordinates of $\abar$ are real and $a_1 + \cdots + a_n \neq 0$ then for any integer $k \leq d/2$, $M_k$, $M_{d-k}$ and $\image \phi_k$ are all equivalent representations.
\begin{proof}
Fix a non-negative integer $k$ for which $k \leq d/2$ and, consequently, $k \leq \lfloor d/2 \rfloor \leq d-k$.
As the kernels of the maps $(\phi_i)_{i \in \mathbb N}$ are decreasing in $i$
and the images of $(\psi_i)_{i \in \mathbb N}$ are increasing, we have
$\ker \phi_{\lfloor d/2 \rfloor} \subseteq \ker \phi_k$ and
$\image \psi_{\lfloor d/2 \rfloor} \subseteq \image \phi_{d-k}$.
So, since $\ker \phi_{\lfloor d/2 \rfloor} + \image \psi_{\lfloor d/2 \rfloor} = V$, we have that
$\ker \phi_{k} + \image \psi_{d-k} = V$.

Finally, $M_k = \image \theta_{F, k}$ and $\theta_{F,k} = \frac{d!}{k!}\,\phi_k \circ \psi_{d-k} \circ \nu_{d-k}$ by Theorem \ref{thm:theta}. Since $\nu_{d-k}$ is surjective, we see that $M_k$ can be expressed more simply as $M_k = \image(\phi_k \circ \psi_{d-k})$. Therefore we have the following chain of equivalences:
\begin{align*}
	M_k
	&= \image(\phi_{k} \circ \psi_{d-k}) \\
	&\cong \image \psi_{d-k} / (\ker \phi_k \cap \image \psi_{d-k}) \\
	&\cong (\ker \phi_{k} + \image \psi_{d-k}) / \ker \phi_k \\
	&\cong V / \ker \phi_{k} \\
	&\cong \image \phi_k.
\end{align*}
As the graded character of $M$ is palindromic, $M_{d-k}$ and $M_k$ are also equivalent.
\end{proof}
\end{prop}

%

\begin{prop} \label{prop:hilbM} Suppose that the coordinates of $\abar$ are real and $a_1 + \cdots + a_n \neq 0$.
The Hilbert function of $M = \langle F \rangle$, the principal submodule of $\apol$ generated by $F \in \apol_d$, is
\[ \hf_{M}(k) = \begin{cases} \rank \phi_k, & {\rm if\ }k \leq d/2;\\
                                             \rank \phi_{d-k}, &{\rm if\ }k \geq d/2
                       \end{cases}
\]
and is unimodal.
\begin{proof}
The degree $k$ homogeneous component of $M$ is
$M_k = \image \theta_k = \image (\phi_k \circ \psi_{d-k})$ by
Theorem \ref{thm:theta} and the fact that $\nu_{d-k}$ is surjective. Thus, $\hf_{M}(k) = \dim_\field M_k = \rank(\phi_k \circ \psi_{d-k})$.
Since $\ker \phi_k$ and $\image \psi_{d-k}$ are transverse,
the rank of the composition of $\phi_k$ and $\psi_{d-k}$ is the minimum of their ranks. By Lemma~\ref{lem:rankpsi}, $\psi_{d-k}$
and $\phi_{d-k}$ have the same rank.   Thus $\hf_{M}(k) =      \min(\rank \phi_k, \rank \phi_{d-k})$.
  Now using Lemma~\ref{lemma:kercontainment}, we see that $\rk \phi_k$ is a non-decreasing function of $k$.  The formula for
  $\hf_M(k)$ now follows.  The unimodality follows immediately.
\end{proof}
\end{prop}

%

\section{The Orbit of a Point and the Algebras $R_\mu$ of DeConcini and Procesi} \label{section:point}
From the results of the previous section, it suffices to know the character and dimension of each $\image \phi_k$
to determine the graded character and Hilbert function of $M$.
In this section, we relate the image of the maps $\phi_k$ to the homogeneous coordinate ring
of the \mbox{$\symn$-orbit} of a projective point. Through this connection, we will express the graded character of $M$ in terms of a Hall-Littlewood Polynomial.

Recall that $L = a_1 x_1 + \cdots + a_n x_n$ is a fixed linear form. Let $p = [a_1: \cdots : a_n] \in \mathbb P^{n-1}$
and $\sigma p = [a_{\sigma^{-1}(1)}: \cdots: a_{\sigma^{-1}(n)}]$ for $\sigma \in \symn$.
The orbit of $p$ is the projective variety $X = \{ \sigma p \mid \sigma \in \symn\}$. Its homogeneous coordinate
ring is $\field[X] = R/I_X$, a one-dimensional arithmetically Cohen-Macaulay ring,
 where $I_X$ is the ideal of homogeneous polynomials vanishing on $X$.

The following elementary lemma gives a sufficient condition for the projective orbit $X$ and the affine orbit $Y = \{ \sigma \abar \mid \sigma \in \symn \}$ to have the same number of points.

\begin{lemma} If $a_1 + \cdots + a_n \neq 0$ then $X$ and $Y$ contain the same number of points.
\begin{proof}
The set $X$ is obtained from $Y$ by identifying affine points which lie on the same line through the origin. Assume that $X$ and $Y$ do not have the same size, so there must be two distinct points $\sigma_i \abar, \sigma_j \abar \in Y$ which represent the same projective point. Thus, for $\tau = \sigma_j^{-1} \sigma_i$,  $\abar$ and $\tau \abar = (a_{\tau^{-1}(1)}, \ldots, a_{\tau^{-1}(n)})$ are distinct points of $Y$, but are equal in $X$. So, there must be some non-zero $z \in \mathbb C$ with $a_i = z a_{\tau^{-1}(i)}$ for all $1 \leq i \leq n$. If $i$ is contained in a cycle of $\tau^{-1}$ of length $m$ then $a_i = z  a_{\tau^{-1}(i)} = z^2 a_{\tau^{-2}(i)} = \cdots = z^m a_i$. Thus, $z$ is an $m$-th root of unity. Also, since $\abar \neq \tau \abar$, we have $z \neq 1$. Therefore, the sum of the coordinates of $\abar$ over the $m$-cycle of $\tau^{-1}$ containing $i$ is
\begin{align*}
	a_i + a_{\tau^{-1}(i)} + \cdots + a_{\tau^{-(m-1)}(i)}
	&= a_i + z^{-1} a_i + \cdots + z^{-(m-1)} a_i\\
	&= a_i (1 + z^{-1}  + \cdots + z^{-(m-1)} )\\
	&= 0.
\end{align*}
Thus, by decomposing $\tau^{-1}$ into its cycles and summing over each cycle, we have expressed $a_1 + \cdots + a_n$ as a collection of disjoint sums which are all zero.
\end{proof}
\end{lemma}

Let $N = \bigoplus_{k \geq 0} N_k \subseteq \apol$ where $N_k = \image \phi_k$. Since $\phi_k:V \to S_k$ is given by $\phi(\sigma_i) = \sigma_i L^k$,
we have $N_k = \spn_\field (\sigma_1 L^k, \ldots, \sigma_\ell L^k)$.

\begin{prop} The annihilator of $N$ is $I_X$ and, furthermore, $N$ and $\field[X] = R/I_X$ are equivalent graded representations of $\symn$.
\begin{proof}
If $\xbar^\bbar$ is a monomial of degree $k$, then $\diff(\xbar^\bbar, L^k) = \bbar! \binom{k}{\bbar} \abar^\bbar = k! \abar^\bbar$ using the multinomial theorem.
Thus, for an arbitrary polynomial $f \in R_k$, we have $\diff(f, L^k) = k! f(\abar)$ by linearity and hence
\[ \diff(f, \sigma_i L^k) = \diff(\sigma_i^{-1}f, L^k) = k! (\sigma_i^{-1}f)(\abar) = k! f(\sigma_i \abar)
\]
for all $1 \leq i \leq \ell$. By the definition of $I_X$, $f \in (I_X)_k$ if and only if $f(\sigma_i \abar) = 0$ for all $1 \leq i \leq \ell$. Thus, $f \in (I_X)_k$ if and only if $\diff(f, \sigma_i L^k) = 0$ for all $1 \leq i \leq \ell$ and hence $(I_X)_k = N_k^\perp = \ann(N)_k$ for all $k \geq 0$. Thus $I_X = \ann(N)$ as both are homogeneous ideals of $R$.

As mentioned in the introduction, the $\symn$-invariant perfect pairings $\diff_k: R_k \times S_k \to \field$ induce equivalences between
 each $(R/\ann(N))_k$ and the dual representation of each $N_k$. Since the $N_k$ are self-dual, $R/I_X$ and $N$ are equivalent graded representations.
\end{proof}
\end{prop}

\begin{prop} If $a_1+\cdots +a_n \neq 0$ then $e_1 = \overline{x_1}+ \cdots +\overline{x_n}$ is not a zero divisor in $\field[X] = R/I_X$. Consequently,
$\field[X]_k \cong \bigoplus_{i=0}^k (\field[X]/(e_1))_i$ as representations for all $k \geq 0$.
\begin{proof}
If $e_1 f \in I_X$ for some $f \in R$, then $e_1(\sigma \abar) f(\sigma \abar) = 0$ for all $\sigma \in \symn$. As $e_1(\sigma \abar) = e_1(\abar) = a_1 + \cdots + a_n \neq 0$, we see that $f(\sigma \abar) = 0$ for all $\sigma \in \symn$. That is, $f \in I_X$. Consequently, $e_1$ is not a zero divisor of $\field[X]$.
Hence, $\field[X]_k$ and $(e_1 \field[X]_{k-1}) \oplus (\field[X]/(e_1))_k$ are equivalent representations. Also,
$e_1 \field[X]_{k-1}$ and $\field[X]_{k-1}$ are equivalent. Thus the result follows by induction.
\end{proof}
\end{prop}

\begin{remark}
Again, let $Y = \{\sigma \abar \in \mathbb A^n  \mid \sigma \in \symn\}$ be the affine orbit of $\abar = (a_1,\ldots, a_n)$ and let $\field[Y] = R/I_Y$
be its (inhomogeneous) coordinate ring. We define the \emph{associated graded algebra} of $\field[Y]$ to be
\[  \gr(\field[Y]) = \bigoplus_{k \geq 0} \field[Y]_{\leq k} / \field[Y]_{\leq k-1}.
\]
Take a non-zero degree $k$ polynomial $f$ and express it as $f = f_k + \cdots + f_0$ where each $f_i$ is homogeneous of degree $i$ and $f_k \neq 0$.
The \emph{leading form} of $f$ is $\LF(f) = f_k$.
The associated graded algebra $\gr(\field[Y])$ is isomorphic to $R/\gr(I_Y)$ where $\gr(I_Y) = \{ \LF(f) \mid f \in I_Y\}$.
One can see this as follows: two degree $k$ polynomials $f,g \in \field[Y]_{\leq k}$ are equal modulo $\field[Y]_{\leq k-1}$ if and only if the leading form of their difference
is in $\gr(I_Y)_k$. Thus, $\field[Y]_{\leq k} / \field[Y]_{\leq k-1}$ is isomorphic to $R_k/\gr(I_Y)_k$. One can also check that this isomorphism is equivariant.
\end{remark}

%

\begin{prop} If $a_1 +\cdots +a_n = t\neq 0$ then $\field[X]/(e_1-t) \cong \field[Y]$.  
Furthermore $\field[X]/(e_1)$ and $\gr(\field[Y])$ are isomorphic graded algebras and equivalent representations of $\symn$.
\begin{proof}
For the first assertion it suffices to show that $I_X + (e_1-t)  = I_Y$.  Since this implies that $I_X + (e_1) = \gr(I_Y)$ the second assertion will also follow.
 
 Clearly $I_X + (e_1-t) \subseteq I_Y$. For the opposite inclusion, consider a non-zero element $f \in I_Y$.  Write 
 $f = f_k + \cdots + f_0 \in I_Y$ where each $f_i \in R_i$ and $f_k \neq 0$.
  We homogenize $f$ with respect to $e_1(\xbar)/t = (x_1+ \cdots + x_n)/(a_1 + \cdots + a_n)$ to
obtain
\[ f' = f_k + \frac{e_1(\xbar)}{t}f_{k-1} + \frac{e_1(\xbar)^2}{t^2}f_{k-2}+ \cdots + \frac{e_1(\xbar)^k}{t^k} f_0.
\]
As $f$ vanishes on $Y$, we see that the homogeneous polynomial $f'$ vanishes on $X$. Thus, $f' \in I_X$.
Since $e_1(\xbar)^s \equiv t^s \pmod{(e_1(\xbar)-t)}$ we see that 
$f \equiv f' \pmod{(e_1(\xbar)-t})$ and thus $f \in I_X + (e_1-t)$.
%
\end{proof}
\end{prop}

The content of the following proposition appears in other works (cf. \cite[Theorem 4.5]{GH2008}), but we include it here for completeness.
\begin{prop} $\field[Y]_{\leq k}$ and $\gr(\field[Y])_{\leq k}$ are equivalent representations.
\begin{proof}
Consider the following short exact sequence of representations:
\[
	 0 \to \field[Y]_{\leq k-1}
	\to \field[Y]_{\leq k} \to \field[Y]_{\leq k}/\field[Y]_{\leq k-1} \to 0.
\]
Since all short exact sequences of representations of finite groups split, $\field[Y]_{\leq k}$ and $(\field[Y]_{\leq k}/\field[Y]_{\leq k-1}) \oplus \field[Y]_{\leq k-1}$ are equivalent. So, by induction on $k$,
$\field[Y]_{\leq k}$ and $\gr(\field[Y])_{\leq k}$ are equivalent representations.
\end{proof}
\end{prop}

Connecting these equivalencies, we have shown that
\begin{prop}
	 \[N_k \cong \field[X]_k \cong (\field[X]/(e_1))_{\leq k} \cong \gr(\field[Y])_{\leq k} \cong \field[Y]_{\leq k}.\]
	\label{prop:NandRmu}
\end{prop}

\begin{remark}
A remarkable fact, proved by Garsia and Procesi, is that $\gr(\field[Y])$ does not depend on the values of the coordinates of $\abar$, but simply on its associated partition $\mu$ \cite[Remark 3.1]{MR1168926}. In view of that, the symbol
\[ R_\mu = \gr(\field[Y])
\]
is used to denote this algebra and $I_\mu = \gr(I_Y)$ denotes the ideal appearing in its presentation as a quotient of $R$.
\end{remark}

The algebra $R_\mu$ has a number of other descriptions.
First, like any artinian algebra, $R_\mu$ is determined by its socle. The socle of $R_\mu$ is the unique irreducible representation of type $\mu$ which appears in the homogeneous component $R_{n(\mu)}$ of degree
\[ n(\mu) = \mu_2 + 2 \mu_3 + \cdots + (r-1) \mu_r.
\]
In fact, $R_{n(\mu)}$ is the lowest degree component of $R$ in which this irreducible representation occurs.

Originally, DeConcini and Procesi \cite{MR629470} defined the ring $R_\mu$ to be the cohomology ring of the variety of flags fixed by a unipotent matrix of shape $\mu = (\mu_1,\ldots, \mu_r)$. They showed that $R_\mu$ could be presented as a quotient of $\mathbb Q[x_1,\ldots,x_n]$ by a homogeneous ideal $I_\mu$ and conjectured a set of generators for $I_\mu$. Tanisaki \cite{MR685425} conjectured a simpler set of generators for $I_\mu$ and, eventually, Weyman \cite{MR1016262} proved these conjectures. Weyman also conjectured a minimal generating set for $I_\mu$, which Biagioli, Faridi and Rosas found to be minimal in some cases and redundant in others \cite{MR2448086}. Garsia and Procesi used Tanisaki's description of $I_\mu$ to show that $R_\mu = \gr(\field[Y])$ as previously mentioned.

We refer the reader to the introduction of \cite{MR1168926} for the progression of papers that led to the graded character of $R_\mu$.
As an ungraded representation, $R_\mu$ is equivalent to the representation afforded by the left cosets of the Young subgroup indexed by $\mu$ or, equivalently \cite[\S5.4]{MR2538310}, to the subrepresentation of $R_{n(\mu)}$ which is spanned by monomials of the form $\prod_{j=1}^r (x_{i_{j,1}} x_{i_{j,2}} \cdots x_{i_{j,\,\mu_j}})^{(j-1)}$ for distinct indices $i_{j,k} \in \{1,\ldots, n\}$. The graded character of $R_\mu$ is given by the combinatorial Hall-Littlewood polynomial:
\begin{equation}
\bigchi_{R_\mu}(t) = Q_\mu(x_x,\dots,x_n;t) = \sum_{\lambda \vdash n} K_{\lambda, \mu}(1/t)t^{n(\mu)} \chi^\lambda
\label{eq:Rmuchar}
\end{equation}
where $K_{\lambda, \mu}(t) \in \mathbb N[t]$ are the Kostka-Foulkes polynomials \cite[Chapter III.6]{MR553598}.
Here $$Q_\mu(x_x,\dots,x_n;t) = t^{n(\nu)}Q'_\mu(x_x,\dots,x_n;t^{-1})$$
where $Q'_\mu(x_1,\ldots,x_n; t)$ is the modified Hall-Littlewood polynomial (see \cite[\S3]{MR1399504}).

The twist in formula (\ref{eq:Rmuchar}) arising from the term $K_{\lambda, \mu}(1/t)t^{n(\mu)}$ where one might expect $K_{\lambda, \mu}(t)$ makes the coefficient of $K_{\lambda, \mu}(t)$ in degree $k$ count the multiplicity of
the irreducible representation of type $\lambda$ that occurs in $(R_\mu)_{n(\mu) - k}$. That is, exponents on $t$ in $K_{\lambda, \mu}(t)$ measure degrees down from the socle of $R_\mu$, rather than up from the constants.

Since $N_k \cong \gr(\field[Y])_{\leq k} \cong (R_\mu)_{\leq k}$, the graded character of $N$ is
\begin{align*}
	\bigchi_N(t)
	&= \sum_{k \geq 0} \bigchi_{(R_\mu)_{\leq k}} t^k\\
	&= \frac{1}{1-t} \sum_{k \geq 0} \bigchi_{(R_\mu)_k} t^k\\
	&= \frac{1}{1-t} \bigchi_{R_\mu}(t).
\end{align*}

Thus we have proved
\begin{thm}\label{thm:maintheorem}
Suppose the coordinates of $\abar$ are real and $a_1 + \cdots + a_n \neq 0$.
 The graded character $\bigchi_M(t) = \sum_{k=0}^d \bigchi_{M_k}t^k$ of $M$
satisfies and is determined by the two equations
$$\bigchi_{M_k} = \bigchi_{M_{d-k}}\quad  {\rm for\ }k \le d/2$$ and
\[ \bigchi_M(t) \equiv \frac{1}{1-t} \bigchi_{R_\mu}(t) \equiv \frac{1}{1-t}Q_\mu(x_1,\dots,x_n;t)   \pmod{ (t^{\lfloor d/2 \rfloor + 1})}.
\]
In particular, $R/I_F$ and $M$ have the same graded character.
\end{thm}

\begin{remark}
 The Hilbert function of $M$ always plateaus at the value of $\ell$. Since $R_\mu = \gr(\field[Y])$ has $\field$-dimension $\ell$ and its socle is in degree $n(\mu)$, by Proposition~(\ref{prop:NandRmu}) we have that $\dim_\field N_k = \dim_\field (R_\mu)_{\leq k} = \ell$ for any $k \geq n(\mu)$. Thus, by Proposition \ref{prop:hilbM}, $\dim_\field M_k = \ell$ for all $k$ with $n(\mu) \leq k \leq d - n(\mu)$. In particular, we need $d \geq 2 n(\mu)$ for one of the homogeneous components of $M$ to achieve dimension $\ell$.

When $d \geq 2 n(\mu)$ it is impossible to write $F$ as a sum of fewer than $\ell$ powers of linear forms since this would force the dimensions of the homogeneous components of $M$ to be smaller than $\ell$. Thus, in this case, the Waring rank of $F$ is $\ell$ and $F = \sum_{i=1}^\ell \sigma_i L^d$ is a Waring decomposition of $F$.
\end{remark}

\begin{remark}
One final observation is that Proposition~(\ref{prop:NandRmu}) implies that the Hilbert function of $N$ is strictly increasing until degree $k = n(\mu)$ and is constant thereafter. Thus, the Hilbert function of $M$ must be strictly unimodal in that it begins with a strictly increasing region, followed by a constant region and ends with a strictly decreasing region (before becoming constant at zero).
\end{remark}

\section{Further Work}

%

The most broad question we propose is to determine the graded characters and Hilbert functions of level artinian quotients of $R$ by $\symn$-stable homogeneous ideals. A graded artinian algebra is level if its socle is contained in a single degree. The ring $R_\mu$ and the coinvariant algebra $R_{1^n}$, in particular, are the most well-studied algebras of this type.
In the introduction, we mentioned the contributions of Bergeron, Garsia and Tesler \cite{MR1779775}, Roth \cite{MR2134302}, and
Morita, Wachi and Watanabe \cite{MR2542134} to this problem.



In this paper, we determined the graded characters of Gorenstein algebras whose socles were spanned by a single  symmetric polynomial $F$ that is the sum of the $\symn$-orbit of a power of a linear form (whose coefficients are real and do not sum to zero).
It remains open to determine the graded character of $R/\ann(F)$ when $F$ is an arbitrary symmetric polynomial.

As mentioned in the introduction, every homogeneous symmetric polynomial $F \in \apol_d$ can be written as a linear combination of orbit sums of powers of linear forms $F_1,\ldots, F_m$ with $F_i = \sigma_1 L_i^d + \cdots + \sigma_{\ell_i} L_i^d$.
We suggest that the graded character of $R/\ann(F)$ may depend on the graded characters of $R/\ann(F_1),\ldots, R/\ann(F_m)$.
If the linear forms determining each $F_i$ are chosen generically and $d$ is sufficiently large, we expect that the character of $R/\ann(F)$ will be the sum of the characters of $R/\ann(F_1), \ldots, R/\ann(F_m)$ in degrees where this is possible.




\bibliographystyle{amsplain}
\bibliography{bibliography}

%
%


\end{document}